# The Friendship Theorem and Minimax Theorems

Dhananjay P. Mehendale
S. P. College, Tilak Road, Pune-411030, India.

## Abstract

We propose a generalization of Hall's marriage theorem. The generalization given here provides a necessary-sufficient condition for arranging a successful friendship among *n* number of *k*-sets. We define multimatrix, multideterminant, multideterminantal monomials, and show that the existence of a nonzero multideterminantal monomial implies the desired existence of a successful friendship. Every nonzero multideterminantal monomial represents a possible way to achieve a successful friendship. A generalization of Menger's theorem and that of some minimax theorems is stated.

**1. Introduction:** The well known theorem of Philip Hall [1] provides a necessary-sufficient condition to achieve a successful marriage of the girls with the boys of their choice as is given in [2]:

**Theorem 1.1**(Hall's theorem-marriage form [2], p. 27)**:** A set of *k* girls can all choose a husband each from the boys they know if and only if any subset of 'r' girls know together (and so can choose one among them) at least 'r' boys.

**Theorem 1.2** (Hall's theorem-matrix form [2], p. 31)**:** Let M be an m×n matrix of zeroes and ones. Then there exists a one in each row of M, no two of which are in the same column, if and only if any set of rows of M (r of them say) have between them ones in at least r columns.

In other words, if
G = { $g_i / 1 \leq i \leq k$ } be the set of girls and let
B = { $b_i / 1 \leq i \leq l$ } be the set of boys and $k \leq l$.
We construct a matrix of size *k* by *l*, as
$A = [a_{ij}]$, where
$a_{ij} = 1$ when girl $g_i$ knows (and so can choose for marrying with) boy $b_j$, and
$a_{ij} = 0$ otherwise.

A successful marriage of all the girls, with the fulfillment of the condition that each girl marrying with the boy of her choice, is thus possible if and only if there exists at least one nonzero determinantal monomial in the determinantal expansion of some $k \times k$ submatrix of $A$. In fact, every such nonzero determinantal monomial represents a successful marriage.

For the sake of simplicity let $l = k$, then the theorem states that for a successful marriage the determinant of the matrix $A$ should contain at least one nonzero determinantal monomial. Note that in a successful marriage we are producing from the (given) $k$-sets (two in number) the 2-sets ($k$ in number) and each 2-set representing a matched pair. We use hereafter an alternative notation for matrices which is more suitable later, Thus, we write $A$ in alternative notation as $<1,2,\cdots,k|1,2,\cdots,k>$)

Suppose there are $n$ number of $k$-sets (of persons), and suppose we want to form $k$ number of $n$-sets, such that exactly one person is chosen from each of the given $k$-sets and all chosen persons are related with each other by (prespecified) friendship bond (or having some desired skill for working together). When is this possible to achieve? We show that this is possible if and only if for every $r$-subset of any set among the $k$-sets there are together at least $r$ elements of every other $k$-set related by the (predefined) friendship bond i.e. for every $r$-subset of any set among the $k$-sets there exist at least one $r$-subset of every other $k$-set having the elements related by the (predefined) friendship bond.

## 2. The Friendship Theorem: Let

$$A^1 = \{a_1^1, a_2^1, \cdots, a_k^1\}$$
$$A^2 = \{a_1^2, a_2^2, \cdots, a_k^2\}$$
$$\vdots$$
$$A^n = \{a_1^n, a_2^n, \cdots, a_k^n\}$$

be $n$ number of $k$-sets.

**Question:** When is it possible to form $k$ number of $n$-sets (we call these sets **friendship sets**) say

$$B^1 = \{a_{j_1}^1, a_{j_1}^2, \cdots, a_{j_1}^n\}$$
$$B^2 = \{a_{j_2}^1, a_{j_2}^2, \cdots, a_{j_2}^n\}$$
$$\vdots$$

$$B^k = \{a^1_{j_k}, a^2_{j_k}, \cdots, a^n_{j_k}\},$$

such that, an element $a^p_{j_q} \in A^p$ and also $a^p_{j_q} \in B^q$, $1 \leq p \leq n$, and $1 \leq q \leq k$, and where all the elements in every $B^q$, $1 \leq q \leq k$, are related by (predefined) friendship bond?

**Theorem 2.1** (The Friendship Theorem): The question raised above has an affirmative answer (i.e. we can form the desired friendship sets) if and only if for any $m$ persons, $1 \leq m \leq k$, in any above given set $A^i$, $1 \leq i \leq n$, there exist together at least $m$ number of friends in every other set $A^j$, $j \neq i$.

Few definitions are now in order.

**Definition 2.1:** A **multimatrix** $Z$ is an $n$ dimensional lattice-like structure, having $k$ elements on each lattice axis, containing in all $k^n$ lattice points, We represent $Z$ as $\langle 1,2,\cdots,k \mid 1,2,\cdots,k \mid \cdots \mid 1,2,\cdots,k \rangle$.

The elements of $Z$ are like: $a_{i_1 i_2 i_3 \cdots i_n}$, such that $1 \leq i_j \leq k$.

**Definition 2.2:** The **multideterminant,** $Det(Z)$, associated with the multimatrix $Z$, is defined as follows:

$$Det(Z) = \sum_{\lambda_1, \lambda_2, \cdots, \lambda_{(n-1)}} \prod_{i=1}^{(n-1)} \mathrm{sgn}(\lambda_i) \prod_{j=1}^{k} a_{j \lambda_1(j) \cdots \lambda_{(n-1)}(j)}$$

where $\lambda_1, \lambda_2, \cdots, \lambda_{(n-1)} \in S_k$, the group of permutations on $k$ symbols, and $\mathrm{sgn}(\lambda_i)$ is the signature of the permutation $\lambda_i$.

**Definition 2.3:** A multideterminantal monomial is a monomial in the multideterminantal expansion given above.

As we have associated a matrix $A$ above (Section 1) with the set of girls $G$, and the set of boys $B$, we now proceed to associate (naturally) a multimatrix with sets: $A^1$, $A^2$, ..., $A^n$. The elements of the multimatrix $Z$ associated with these sets $a_{i_1 i_2 i_3 \cdots i_n}$ have two values:

$a_{i_1 i_2 i_3 \cdots i_n} = 1$, if in the $n$-set of elements $\{a^1_{i_1}, a^2_{i_2}, \cdots, a^n_{i_n}\}$ where $a^j_{i_j} \in A^j$ for all $j$ and where all elements $a^j_{i_j}$, $1 \leq j \leq n$, are related by (predefined) friendship bond, i.e. set $\{a^1_{i_1}, a^2_{i_2}, \cdots, a^n_{i_n}\}$ forms a set like, $B^i, i = 1, 2, \cdots, k$, and $a_{i_1 i_2 i_3 \cdots i_n} = 0$, otherwise.

Note that there are in all $(k!)^{(n-1)}$ multideterminantal monomials in the multideterminant of the multimatrix having $k^n$ elements. Each nonzero element of the multimatrix represents a set like $B^i, i = 1, 2, \cdots, k$. Each nonzero multideterminantal monomial gives all the desired friendship sets $B^i, i = 1, 2, \cdots, k$.

**Lemma 2.1:** The question raised above has an affirmative answer if and only if in the associated multimatrix there exists at least one nonzero multideterminantal monomial.

**Proof:** When a nonzero monomial exists we can form the desired sets $B^i, i = 1, 2, \cdots, k$, and when the desired sets exist we can form the nonzero monomial. We are essentially using the following association, namely, $a_{i_1 i_2 i_3 \cdots i_n} = 1 \Leftrightarrow \{a^1_{i_1}, a^2_{i_2}, \cdots, a^n_{i_n}\}$ is a set of friends, called a friendship set. □

**Proof of theorem 2.1:** Necessity is straightforward since if any member of any set $A^i$ finds a friend in every other set $A^j, j \neq i$, then clearly any $m$ members of any set $A^i$ must find at least $m$ friends in every other set $A^j, j \neq i$.

Suppose for any $m$ persons, $1 \leq m \leq k$, in any above given set $A^i$, $1 \leq i \leq n$, there exist together at least $m$ number of friends in every other set $A^j$, $j \neq i$. We proceed to show that a nonzero multideterminantal monomial can be formed. Suppose not, i.e. there doesn't exist a nonzero multideterminantal monomial. It meant that if we try to form a nonzero multideterminantal monomial (in all possible ways), by selecting in succession the elements of the multimatrix $a_{i_1 i_2 i_3 \cdots i_n}$, such that $a_{i_1 i_2 i_3 \cdots i_n} = 1$, then process halts at some step $m < k$. Hence, any process of selection of elements of the multimatrix like $a_{i_1 i_2 i_3 \cdots i_n} = 1$, in order to construct a nonzero multideterminantal monomial, terminates at some step $m < k$. It essentially implies that the process of forming a nonzero multideterminantal monomial halts due to arrival of a multideterminant of the multisubmatrix of size $m$, containing in all $m^n$ lattice points, such that every multideterminantal monomial contained therein has zero value. Since we reach to the same outcome for all choices of constructing a nonzero multideterminantal monomial this in turn implies that for some $m$ persons, $1 \leq m \leq k$, in a set $A^i$, $1 \leq i \leq n$, there exist together fewer than $m$ number of friends in some other set $A^j$, $j \neq i$, a contradiction.

Once a nonzero multideterminantal monomial is formed the rest is clear from the lemma 2.1. Hence, the theorem. □

**3. The Graph Theoretic Interpretation:** This result can be interpreted in the language of graph theory as follows:

Let G be a graph with vertex set $V(G) = \bigcup_{i=1}^{n} A^i$, and the the friendship relation between any two points belonging to $A^i$ and $A^j$, $j \neq i$ for all $i$ and $j$ be the edges of G. Then it is easy to see that the theorem offers a condition for **decomposition** of **this** graph G into **disjoint cliques** on n points, i.e. such a graph can be **factored** into k number of disjoint cliques on n points.

**4. Minimax Theorems:** In this section as an application of **ubiquitous Hall's theorem** we discuss the theorem of Konig-Egervary. In its matrix form it concerns taking a matrix of $0^s$ and $1^s$ and looking for one 1 in each row with no two in the same column. However, such a set of $1^s$ may not exist and we may instead look for as many $1^s$ as possible with no two in the same column. As an example, we see here the generalization of the well known minimax theorem, the theorem of Konig-Egervary (matrix and graph form). In a matrix there are rows and columns and every element $a_{ij}$ belongs to $i^{th}$ row and $j^{th}$ column. The rows are along one dimension while the columns are along the other dimension when one views a matrix as a two dimensional object. We call hereafter a row or a column **a line**. In case of an n dimensional multimatrix every element (entry) $a_{i_1 i_2 i_3 \cdots i_n}$ belongs to n lines $\{i_1, i_2, \cdots, i_n\}$. The usual matrix of size $k$ (i.e. $k \times k$) contains together $2k$ lines ($k$ rows and $k$ columns) and can be denoted as (is done above) $<1,2,\cdots,k | 1,2,\cdots,k>$. The n dimensional multimatrix of size $k$ (i.e. $k \times k \times \ldots \times k$, $k$ taken n times) contains $nk$ lines and denoted as (is done above)

$\langle 1, 2, \cdots, k \mid 1, 2, \cdots, k \mid \cdots \mid 1, 2, \cdots, k \rangle$. We now give the following generalization of the theorem of **Konig-Egervary** [2].

**Theorem 4.1:** Let M be an n dimensional multimatrix of size $k$ of 0s and 1s. Then the minimum number of lines containing all the 1s of M is equal to the maximum number of 1s with no two in the same line.

**Proof:** Let $\alpha$ be minimum number of lines of M which between them contain all its 1s and let $\beta$ be the maximum number of 1s of M with no two in the same line.

Since there are $\beta$ 1s with no two in the same line it takes at least $\beta$ lines to include those 1s, hence $\alpha \geq \beta$.

To show that $\alpha \leq \beta$, we shall find $\alpha$ 1s with no two in the same line. Let $\alpha = l_1 + l_2 + \cdots + l_n$, where without loss of generality (since we can manage this by appropriate **line interchanges** if it is not so) the first $l_i$ lines along the $i^{th}$ dimension of M, $i = 1, 2, \ldots, n$; among them contain all the 1s of M. Consider the multimatrix $N_1$ made up of **first** $l_1$ lines along $1^{st}$ dimension and (k - $l_2$), (k - $l_3$), …, (k - $l_n$) **later** lines (i.e. lines after first $l_i$ lines along the $i^{th}$ dimension of M, $i = 2, 3, \ldots, n$) along $2^{nd}, \ldots, n^{th}$ dimension respectively. We claim that any r lines among the **first** $l_1$ lines (along $1^{st}$ dimension) have between them 1s in at least r lines among (k - $l_2$), (k - $l_3$), …, (k - $l_n$) **later** lines along $2^{nd}, \ldots, n^{th}$ dimension respectively. For if not there would be some r lines among the **first** $l_1$ lines along $1^{st}$ dimension with 1s in just u lines along other dimensions, where u < r. But then replacing the r lines of M along $1^{st}$ dimension by these u lines along other dimensions would give a smaller set of lines containing all the 1s of M, which is not possible. Hence any r lines of N along the $1^{st}$ dimension do contain 1s in at least r lines along other dimensions. Therefore N has a set of $l_1$ 1s with one in each line along $1^{st}$ dimension and no two in the same line along any other dimension. Continuing with the **similar argument** we Consider the multimatrix $N_2$ made up of **first** $l_2$ lines along $2^{nd}$ dimension and (k - $l_1$), (k - $l_3$), …, (k - $l_n$) **later** lines along $1^{st}, 3^{rd}, \ldots, n^{th}$ dimension respectively and we show in the same way that any r lines among the **first** $l_2$ lines (along $2^{nd}$ dimension) have between them 1s in at least r lines among (k - $l_1$), (k - $l_3$), …, (k - $l_n$) **later** lines along $1^{st}, 3^{rd} \ldots, n^{th}$ dimension respectively. Thus, we can see that $N_2$ has a set of $l_2$ 1s with one in each line along $2^{nd}$ dimension and no two in the same line along any other dimension. Continuing with this argument we see that there are

$\alpha$ 1s in M with no two in the same line. Since $\beta$ is the maximum number of such 1s, it follows that $\alpha \le \beta$.  $\square$

We have defined above the **lines** of a multimatrix. In precise language, if we fix any $(k-1)$ indices and allow one index of an element of a multimatrix to vary we get set of multimatrix elements belonging to a **line** of the multimatrix. If we fix any $(k-2)$ indices and allow two indices of an element of a multimatrix to vary we get the set of multimatrix elements belonging to a (coordinate) **2-plane** of the multimatrix. Likewise, if we fix any $(k-r)$ indices, $r < k$, and allow $r$ indices of an element of a multimatrix to vary we get the set of multimatrix elements belonging to a (coordinate) ***r*-plane** of the multimatrix.

We now state the following generalization of theorem 4.1(which can be settled by proceeding as per the proof of theorem 4.1):

**Theorem 4.2:** Let M be an n dimensional multimatrix of size $k$ of 0s and 1s. Then the minimum number of (distinct) *r*- planes containing all the 1s of M is equal to the maximum number of 1s with no two are in the same *r*-plane.

**Remark 4.1:** Theorem 4.1 can be used for generalizing the well known Hungarian Method [3], [4] for multimatrices to solve problems in Operations Research of the following type:

**Minimize:** $Z = \sum_{i_1, i_2, \cdots, i_k} C_{i_1 i_2 \cdots i_r} x_{i_1 i_2 \cdots i_k}$

**Subject to:** $\sum_{i_r} x_{i_1 i_2 \cdots i_k} = 1,$ for all $r = 1, 2, \ldots, k$

**Remark 4.2:** The following theorem is virtually a **restatement** of the above given generalized version of the theorem of Konig-Egervary:

**Theorem 4.3:** Let G = ($V_1, V_2, \cdots, V_k, E$) be a multipartite graph. Then the minimum number of vertices which between them include at least one end point of each edge is equal to the maximum number of edges in a matching.

**Definition 4.1:** A path from any vertex subset $V_i$ to $V_j$, for all $i$ and $j$, $i \ne j$ is simply a path whose first vertex is in $V_i$ and whose last vertex is in $V_j$. Paths with no vertex in common are called disjoint paths.

**Definition 4.2:** Two paths, one from any vertex subset $V_i$ to $V_j$ and the other from any vertex subset $V_k$ to $V_l$, $i \neq j$ and $k \neq l$ and $i \neq k$ and/or $j \neq l$ are said to be disjoint if they have no vertex in common.

**5. Menger's Theorem:** Perhaps the **most important theorem** in the entire graph theory is Menger's theorem, the theorem first proved by K. Menger in 1927. One can obtain the following generalization of this theorem by proceeding on similar lines as is done in [2], [4] for the case $k = 2$. There are many important theorems that directly follow from the Menger's theorem and its corollaries. The excellent account of this development can be found in [4].

We now conclude our discussion with a statement that the theorem 4.3 is true for **any** $k$ subsets of the vertex set in **any** graph as is done by Menger for the case $k = 2$.

**Theorem 5.1(Generalized Menger's Theorem):** Let $G = (V, E)$ be any graph. Let $V_1, V_2, \cdots, V_k$ be disjoint subsets in the partitioning of $V$. Then, the minimum number of vertices taken together that separate all subsets $V_i$ from $V_j$, for all $i$ and $j$, where $i \neq j$ equals the maximum number of disjoint paths taken together from $V_i$ to $V_j$, for all $i$ and $j$, $i \neq j$.

## Acknowledgements

I am thankful to Dr. M. R. Modak, Dr. S. A. Katre, Bhaskaracharya Pratishthana, Pune, and Dr. T. T. Raghunathan, University of Pune, for their keen interest and useful discussions.